\documentclass[10pt]{article}

\catcode `\@=11
\@addtoreset{equation}{section}

\begin{document}

\hsize 140mm
\textheight 190mm

\title{\bf Groups with involution, and quasigroups with cracovian representations} 
\author{Jerzy Koci\'nski}

\maketitle

\footnote{Professor J. Koci\'nski, Grzybowska 5, m. 401, 00-132 Warszawa, Poland\\
e-mail: kocinsk@if.pw.edu.pl}

\newpage

\noindent
{\bf Abstract}$\,$ In groups with involution a nonassociative product of elements is defined, which leads to the definition of a certain type of quasigroups. These quasigroups are represented by square tables of complex numbers, with inverses, which differ from the matrix representations of groups in the rule of performing the product of two tables. The row-by-column product of two matrices in representations of groups is replaced by the column-by-column product, which is called the cracovian product, in representations of the defined type of quasigroups. The matrices undergoing the column-by-column product are called cracovians. The basic properties of the quasigroups connected with groups with involution are determined while only a summary of the properties of cracovian algebra is presented, as the basis of cracovian representation theory for the quasigroups connected with groups with involution. Clifford groups are groups with involution and the quasigroups connected with them are determined. The orthogonal and pseudo-orthogonal rotation groups belong to groups with involution.
An analogy is drawn between Weyl's "hidden" symmetry group of an object, and the quasigroup connected with the group with involution.

\vspace{7mm}

\section{\bf Introduction}

\noindent
The work in the field of quasigroups acquired impetus with the investigations presented in the "Projektive Ebenen" by G. Pickert (1955)\cite{Pickert}, in the "Binary Systems" by R. H. Bruck (1958 and 1966)\cite{Bruck}, and in the "Foundations of Quasigroups and Loops" 
by V. D. Belousov (1967) \cite{Belousov}. It was followed by the books which appeared in the last decade of the twentieth century: "Quasigroups and Loops: Introduction" 
by H. O. Pflugfelder (1990)\cite{Pflugfelder}, "Quasigroups and Loops: Theory and Application" edited by O.Chein, H. O. Pflugfelder and J. D. H. Smith (1990) \cite{Chein}, and "Smooth Quasigroups and Loops" by L. V. Sabinin (1999)\cite{Sabinin1}.  A survey of papers on nonassociative geometry with the reference to space-time was presented by L. V. Sabinin \cite{Sabinin2}. The line of thought of that survey was pursued in \cite{Sbitneva}, in an application of nonassociative geometry to special relativity. Quasigroups connected with groups with involution, and their cracovian representations were introduced in  \cite{Kocinski1, Kocinski2, Kocinski3}. The attempts of applying nonassociative algebras in quantum mechanics, were initiated by Jordan \cite{Jordan1}, and Jordan {\em et al.,} \cite{Jordan2, Jordan3}.  Full expositions of Jordan algebras were given by Braun and Koecher \cite{Braun} and by Jacobson \cite{Jacobson}.
The nonassociative and noncommutative column-by-column product, (or, alternatively, row-by-row product) of two tables of symbols or complex numbers, was defined and applied in
theoretical astronomy by Banachiewicz in the year 1924, at the Jagellonian University in Cracow. 
The name "cracovian"  was coined in connection with that. The cracovian algebra was developed by Banachiewicz in a series of papers over the span of thirty years, of which we quote the papers \cite{Banachiewicz1, Banachiewicz2, Banachiewicz3, Banachiewicz4, Banachiewicz5, Banachiewicz6, Banachiewicz7, Banachiewicz8, Banachiewicz9, Banachiewicz10}.
It was presented in the books by Banachiewicz \cite{Banachiewicz11}, by Sierpi\'nski \cite{Sierpinski}, by Lukaszewicz and Warmus \cite{Lukaszewicz}, and by Koci\'nski \cite{Kocinski3}. In the last book a list of references concerning cracovian algebra and some of its applications are given.

\section{The quasigroups connected with groups with involution}

\noindent
We consider a finite group $G$ of order $N$, with elements $g$ and unit element $e$, for which an involution
operation $I(g)$ exists with the following properties:
\begin{eqnarray}
I(g)\in G, \quad\forall g\in G
\nonumber\\
I[I(g)]=g
\nonumber\\
I(g_1 g_2 ... g_p)=I(g_p)I(g_{(p-1)})...I(g_2)I(g_1), \quad\forall g_1, g_2, ..., g_p\in G
\nonumber\\
I(\pm e)=\pm e 
\label{eq2:1}
\end{eqnarray}

\noindent
In the following we will denote the group elements by 
$g_1, g_2,...,g_N$ or, alternatively, by $a, b, c,...,z$.
\\  
\noindent
{\bf Definition 2.1.}  We define the following group automorphism by the equality

\begin{equation}
I(a)b\buildrel\rm def\over=c, \quad \mbox{for a fixed}\quad a,\quad\forall a,b,c\in G
\label{eq2:2}
\end{equation}

\noindent
where, for brevity, $a,b$ and $c$ now denote arbitrary products of the elements $g_1, g_2,...,g_N$.
\\
\noindent
{\bf Corollary 2.1.}  The automorphisms defined in Eq. (\ref{eq2:2}) together with
Eq. (\ref{eq2:1}) form groups.
\\
\noindent
Proof. The product of two such automorphisms yields another automorphism of the same type,
$$
I(d)[I(a)b]=I(d)I(a)b=I(ad)b
$$

\noindent
The product of three automorphisms is associative since we have
$$
[I(c)I(b)]I(a)d=I(bc)I(a)d=I(abc)d
$$
$$
I(c)[I(b)I(a)d]=I(c)I(ab)d=I(abc)d
$$

\noindent
There exists the unit element since $I(1)a=a=aI(1)$,
and for each automorphism $I(a)$, there exists the inverse automorphism
$I(a^{-1})$ since we have
$$
I(a^{-1})[I(a)b]=I(aa^{-1})b=b
$$
$$
I(a)[I(a^{-1})b]=I(a^{-1}a)b=b
$$

\noindent
{\bf Definition 2.2}.  The "dot" product denoted by $(\cdot)$ is defined by

\begin{equation}
a\cdot b\buildrel\rm def\over=I(b)a,\qquad\forall a,b \in G
\label{eq2:3}
\end{equation}

\noindent
where on the right hand side we are dealing with the associative product in the group $G$.
\\
\noindent
{\bf Observation 2.1.}
To distinguish the elements of the group $G$ undergoing the associative product from the same
elements in the "dot' product,
we will write the elements of the group $G$ in the "dot" product with a caret. This means that
the product $a\cdot b$ will be rewritten as $\hat a\cdot\hat b$. 
\\
\noindent
{\bf Definition 2.3.}  If $e$ denotes the unit element in the group $G$ in
Eq. (2.1), and $a$ is an element
in that group, the right unit element $\tau$ for the "dot" product in
Eq.(2.3), is defined by the equality

\begin{equation}
\hat a\cdot \tau=I(e)a=a\rightarrow \hat a, \qquad\forall a\in G
\label{eq2:4}
\end{equation}

\noindent
where the arrow on the right hand side indicates the passage from the group element $a$
to the same element undergoing the "dot" product, then denoted by $\hat a$.
\\
\noindent
{\bf Corollary 2.2.}  The right unit element $\tau$ is not at the same time a left
unit element, since from Eq. (2.3) we obtain

\begin{equation}
\tau\cdot\hat a=I(a)e=I(a)
\label{eq2:5}
\end{equation}

\noindent
and, in a general case, we have $I(a)\neq a$.
\\
\noindent
{\bf Corollary 2.3.}  The right unit element $\tau$ applied from the left twice, does
not change any element in the "dot" product, since we have

\begin{equation}
\tau\cdot (\tau\cdot\hat a)=I[I(a)e]e=a\rightarrow \hat a, \qquad\forall a\in G 
\label{eq2:6}
\end{equation}

\noindent
{\bf Definition 2.4.}  The product defined in Eq. (2.3)
is performed from left to right:

\begin{equation}
\hat a\cdot\hat b\cdot\hat c\cdot ... \cdot\hat z=\{[(\hat a\cdot\hat b)\cdot\hat c]
\cdot ... \cdot\hat z\}
\label{eq2:7}
\end{equation}

\noindent
{\bf Corollary 2.4.}  From Eq. (2.3) it follows that

\begin{equation}
(\hat a\cdot\hat b)\cdot\hat c=\hat a\cdot [\hat c\cdot (\hat\tau\cdot b)], \qquad\forall
a,b,c\in G
\label{eq2:8}
\end{equation}

\noindent
Proof. We have:
$$
(\hat a\cdot\hat b)\cdot\hat c=I(c)(\hat a\cdot\hat b)=I(c)I(b)a
$$
\noindent
and
$$
\hat a\cdot [\hat c\cdot (\hat\tau\cdot\hat b)]=I[I(\hat\tau\cdot\hat b)c]a=I(c)I(b)a
$$

\noindent
The "dot" product is nonassociative.
\\
\noindent
{\bf Corollary 2.5.}  From Eqs. (2.3) and (2.5) it follows that

\begin{equation}
ab=\hat b\cdot I(a), \qquad\forall a,b\in G
\label{eq2:9}
\end{equation}

\noindent
where on the right hand side we have $I(a)=a'\rightarrow \hat a'$.
\\
\noindent
{\bf Definition 2.5.}  If $a^{-1}$ is the inverse of the element $a$ in the group $G$,
then the right inverse of $\hat a$ is defined by

\begin{equation}
\hat a^{-1}=I(a^{-1})
\label{eq2:10}
\end{equation}

\noindent
since from Eq. (2.10) we then obtain: $e=a^{-1}a=
\hat a\cdot I(a^{-1})\rightarrow\tau$.
\\
\noindent
{\bf Corollary 2.6.} The right inverse in Eq. (2.10) at the same time is the left
inverse, since from Eqs. (2.7) and (2.10) we obtain

\begin{equation}
\tau=\tau\cdot\tau=\tau\cdot (\hat a\cdot\hat a^{-1})=
[\tau\cdot (\tau\cdot\hat a^{-1})]\cdot\hat a=\hat a^{-1}\cdot\hat a
\label{eq2:11}
\end{equation}

\noindent
{\bf Corollary 2.7.} The right or left inverse of $\tau$ 
is equal to $\tau$.
\\
\noindent
{\bf Corollary 2.8.} For the "dot" product there is no left unit element.
\\
\noindent
Proof. If $\hat x$ were a left unit element, we would have $\hat x\cdot 
(\hat b\cdot\hat c)=\hat b\cdot\hat c.$ However, from Eq. (2.6) at the same
time we would have $\hat x\cdot (\hat b\cdot\hat c)=
[\hat x\cdot (\tau\cdot\hat c)]\cdot\hat b=
(\tau\cdot\hat c)\cdot\hat b.$ We therefore obtain: $\hat b\cdot\hat c=
(\tau\cdot\hat c)\cdot\hat b$, which, in a general case, is not true.
\\
\noindent
{\bf Corollary 2.9.} From Eqs. (2.2), (2.3),
(2.6) (2.7) and (2.9) it follows that
the elements of the group $G$ in Eq. (2.1) undergoing the "dot" product,
form a quasigroup $QG$ with the right unit element.
The Cayley table of this quasigroup is determined on the basis
of the product definition in Eq. (2.3) and the Cayley table of the group $G$. 
\\
\noindent
{\bf Observation 2.2.} If the "dot" product definition in Eq.(2.2) were replaced
by the definition

\begin{equation}
\hat a\cdot\hat b\buildrel\rm def\over=aI(b), \qquad\forall a,b \in G 
\label{eq2:12}
\end{equation}

\noindent
we again would obtain Eqs. (2.3), (2.6), (2.7) and (2.9), 
which define a quasigroup $QG$, connected with the group $G$. 
The choice between the two definitions of the "dot"
product will be dictated by the homomorphism existing between the quasigroup based on
Eq. (2.3), and square cracovians with inverses, and the lack of such a homomorphism
when the product definition in Eq. (2.12) is accepted.
\\
\noindent
{\bf Observation 2.3}.
In the following formulas the definition in Eq. (2.3) which relates the quasigroup "dot"
product with the group product will not be used. To simplify notation, we no longer will 
write the "caret" above the elements of the group $G$ undergoing the quasigroup product. 
\\
\noindent
{\bf Observation 2.4.} 
In order to avoid unnecessary brackets in the formulas, we will omit from now on the
"dot" while multiplying any element of the quasigroup from the left or from the right
by the right unit element $\tau$. Consequently, instead of $\tau\cdot a$ or $a\cdot\tau$,
we will write $\tau a$ or $a\tau$, respectively; instead of $b\cdot (\tau\cdot a)$ we will
write $b\cdot\tau a$, and $a\cdot\tau\cdot b$ will be replaced by $a\tau\cdot b$. 
In particular, we now will rewrite Eq. (2.7) in the form

\begin{equation}
(a\cdot b)\cdot c=a\cdot (c\cdot\tau b)
\label{eq2:13}
\end{equation}

\noindent
{\bf Corollary 2.10.} From Eqs. (2.13) and (2.6) we find that

\begin{equation} 
\tau (a\cdot b)=b\cdot a
\label{eq2:14}
\end{equation}

\noindent
{\bf Corollary 2.11}. From $a\cdot a^{-1}=\tau$, it follows that 
$(\tau a)\cdot (\tau a^{-1})=\tau$.
\\
\noindent
Proof. Writing $a^{-1}=b$, and 
multiplying $a\cdot b=\tau$ from the right by $(\tau b)^{-1}$ we obtain $(a\cdot b)\cdot
(\tau b)^{-1}=a\cdot [(\tau b)^{-1}\cdot\tau b]=\tau(\tau b)^{-1}=a$, since $a\cdot b=\tau$.
Multiplying the last equality from the left by $\tau b$ we obtain $(\tau b)^{-1}=\tau a$,
and multiplying both sides of the last equality by $\tau b$ from the right we obtain,

\begin{equation}
(\tau b)^{-1}\cdot (\tau b)=(\tau a)\cdot (\tau b)=(\tau a)\cdot (\tau a^{-1})=\tau
\label{eq2:15}
\end{equation}

\noindent
which proves the statement.
\\ 
\noindent
{\bf Corollary 2.12.}  We have

\begin{equation}
(a\cdot b)^{-1}=a^{-1}\cdot b^{-1}
\label{eq2:16}
\end{equation}

\noindent
Proof. From Eqs. (2.13) and (2.6) we obtain

$$
(a\cdot b)\cdot (a^{-1}\cdot b^{-1})=a\cdot [(a^{-1}\cdot b^{-1})\cdot\tau b]=
$$
$$
a\cdot [a^{-1}\cdot (\tau b\cdot\tau b^{-1})]=a\cdot a^{-1}=\tau
$$

\noindent
which proves the statement. 
\\
\noindent
{\bf Corollary 2.13.}  From Eq. (2.16) we obtain

\begin{equation}
(\tau a)^{-1}=\tau a^{-1}
\label{eq2:17}
\end{equation}

\noindent
since $\tau^{-1}=\tau$.
\\
\noindent
{\bf Corollary 2.14.} For any finite number of factors $g_1,g_2,...,g_p$, we have,

\begin{equation}
\Big(g_1\cdot g_2\cdot ... \cdot g_p\Big)^{-1}=g^{-1}_1\cdot g^{-1}_2\cdot ... \cdot
g^{-1}_p
\label{eq2:18}
\end{equation}

\noindent
since if this equality is valid for $k$ factors, then owing to 
Eq. (2.16) it also is valid for $(k+1)$ factors,

$$
[(g_1\cdot g_2\cdot ... \cdot g_k)\cdot g_{(k+1)}]^{-1}=
(g_1\cdot g_2\cdot ... \cdot g_k)^{-1}\cdot g^{-1}_{(k+1)}
$$
$$
g^{-1}_1\cdot g^{-1}_2\cdot ... \cdot g^{-1}\cdot g^{-1}_{(k+1)}
$$

\noindent
{\bf Corollary 2.15.} From Eqs. (2.13) and (2.14) we find that

\begin{equation}
\tau (a\cdot b\cdot c)=c\cdot \tau b\cdot a
\label{eq2:19}
\end{equation}

\noindent
since from Eq. (2.14) we have $\tau [(a\cdot b)\cdot c]=c\cdot (a\cdot b)$, and
from Eq. (2.17) we obtain $c\cdot (a\cdot b)=c\cdot\tau b\cdot a$.
\\
\noindent
{\bf Corollary 2.16.}  For any finite number of $p$ factors we have

\begin{equation}
\tau\Big(g_1\cdot g_2\cdot ... \cdot g_{(p-1)}\cdot g_p\Big)=
g_p\cdot \tau g_{(p-1)}\cdot ... \cdot\tau g_2\cdot g_1
\label{eq2:20}
\end{equation}

\noindent
Proof. We write: $g_1\cdot g_2\cdot ... \cdot g_p=(g_1\cdot g_2\cdot ... \cdot 
g_{(p-2)})\cdot g_{(p-1)}\cdot g_p=W$, and then owing to Eqs. (2.19) and (2.13)
we obtain the equalities

\begin{small}
$$
\tau W=g_p\cdot\tau g_{(p-1)}\cdot (g_1\cdot g_2\cdot ... \cdot g_{(p-2)})=
(g_p\cdot \tau g_{(p-1)})\cdot [(g_1\cdot g_2\cdot ... \cdot g_{(p-3)})\cdot g_{(p-2)}]=
$$
$$
(g_p\cdot\tau g_{(p-1)}\cdot\tau g_{(p-2)})\cdot [(g_1\cdot g_2\cdot ... \cdot g_{(p-4)})
\cdot g_{(p-3)}]=...=g_p\cdot\tau g_{p-1}\cdot ...\cdot\tau g_2\cdot g_1
$$
\end{small}

\noindent
{\bf Clifford quasigroups} \cite{Kocinski1, Kocinski2, Kocinski3}. As an example of quasigroups connected with  groups with involution we will consider the quasigroups connected with Clifford groups.
We consider Clifford algebras with the generators $\gamma_1, \gamma_2,..., \gamma_N,\,
N=1,2,...,$, and with the structural condition:

\begin{equation}
\gamma_{\mu}\gamma_{\nu}+\gamma_{\nu}\gamma_{\mu}=2\delta_{\mu\nu},\qquad
\mu,\nu=1,2,...,N
\label{eq2:22}
\end{equation}

\noindent
For a fixed $N$, the basis of a Clifford algebra consists of the unit element $1$,
the generators $\gamma_1,...,\gamma_N$, and all linearly independent products of these 
generators. The dimension of this algebra is $2^N$. Since the respective Clifford group $G$
contains the element $(-1)$, the order of the group $G$ is equal to $2^{N+1}$.
\\
\noindent
{\bf Definition 2.6.}  The involution operation in a Clifford group $G$ is defined by
\begin{eqnarray}
I(\gamma_{\mu})\buildrel\rm def\over=-\gamma_{\mu}, \quad I(I(\gamma_{\mu}))\buildrel\rm def\over=\gamma_{\mu}, \quad I(\pm 1)\buildrel\rm def\over=\pm 1
\nonumber\\
I(\gamma_{\mu}\gamma_{\nu} ... \gamma_{\sigma})\buildrel\rm def\over=I(\gamma_{\sigma}) ... I(\gamma_{\nu})
I(\gamma_{\mu}), \quad \gamma_{\mu}, \gamma_{\nu}, ... ,\gamma_{\sigma}\in G
\label{eq2:22}
\end{eqnarray}

\noindent
For the group elements consisting of products of $\gamma 's$, we define the notation

\begin{equation}
\gamma_{\mu}\gamma_{\nu} ... \gamma_{\sigma}\buildrel\rm def\over=\gamma_{\mu\nu ...\sigma}
\label{eq2:23}
\end{equation}

\noindent
{\bf Definition 2.7.}  We define the following group automorphism by the equality:

\begin{equation}
I(\gamma_A)\gamma_B\buildrel\rm def\over=\gamma_C, \quad \mbox{for a fixed}\quad \gamma_A
\label{eq2:24}
\end{equation}

\noindent
where, for brevity, $\gamma_A, \gamma_B,$ and $\gamma_C$ now denote arbitrary elements
of a Clifford group, i.e., also arbitrary products of the elements $\gamma_{\alpha},\,
\alpha=1,2,...,N$. 
The inverse of any $\gamma_A=\gamma_{\alpha}\gamma_{\beta}...\gamma_{\rho}$ is equal to
$\gamma^{-1}_A=\gamma_{\rho}...\gamma_{\beta}\gamma_{\alpha}$.
\\
\noindent
{\bf Corollary 2.17.}  The automorphisms defined in Eq. (\ref{eq2:24}) together with
Eq. (\ref{eq2:22}) form groups.
\\
\noindent
{\bf Corollary 2.18}.  The automorphism defined in Eq. (2.25) preserves the 
structural condition in Eq. (2.21).
\\
\noindent
{\bf Definition 2.8.}  We define the product

\begin{equation}
(\hat\gamma_{\mu\nu...\sigma})\cdot (\hat\gamma_{\epsilon\eta...\rho})\buildrel\rm def\over=I(\gamma_{\epsilon
\eta...\rho})\gamma_{\mu\nu...\sigma}
\label{eq2:26}
\end{equation} 

\noindent
where the product on the right hand side is the associative product in Clifford groups,
and where the $\gamma-$symbols in the "dot" product are distinguished with a "caret" from
the same $\gamma-$symbols in the associative product.
For single$-\gamma$ group elements this definition reads

\begin{equation}
\hat\gamma_{\mu}\cdot\hat\gamma_{\nu}\buildrel\rm def\over=I(\gamma_{\nu})\gamma_{\mu}
\label{eq2:27}
\end{equation}

\noindent
Eqs. (2.4) through (2.26) determine Clifford quasigroups connected with Clifford groups.
\\
\noindent
{\bf Observation 2.14.} A respective nonassociative Clifford algebra has a basis consisting
of the right unit element $\tau$, the generators $\hat\gamma_{\mu},\,\mu=1,...,N$,
fulfilling the condition 

\begin{equation}
\hat\gamma_{\mu}\cdot\hat\gamma_{\nu}+\hat\gamma_{\nu}\cdot\hat\gamma_{\mu}=
-2\delta_{\mu\nu}\tau, \quad \mu,\nu=1,...,N
\label{eq2:28}
\end{equation}

\noindent
and of all linearly independent products of these generators. The dimension of this algebra
is $2^N$. To perform the multiplication of these generators or of their products, we have to 
apply Eqs. (2.26) or (2.25). Since the respective Clifford quasigroup contains the element $-\tau$, the order of that quasigroup is equal to $2^{N+1}$.

\section{The Basic Properties of Cracovian Algebra}

\noindent
The proofs of the basic properties of cracovian algebra can be found in \cite{Kocinski3}.
\\
\noindent
{\bf Observation 3.1.} A square or a rectangular table of symbols or complex numbers can be 
called a matrix or a cracovian, depending on the definition of the product of two such 
tables.
The definitions of a product with a scalar, of a sum or a difference
of two matrices, of the equality of two matrices,
of a symmetric or an antisymmetric matrix, and
of a transposed matrix, therefore carry over from matrices on cracovians.
\\
\noindent
{\bf Definition 3.1.} A rectangular cracovian $A$ is defined as the table of $m\times
n$ elements

\begin{equation}
\left\{
\begin{array}{rrrr}
a_{11} & a_{21} & ... & a_{m1} 
\\
a_{12} & a_{22} & ... & a_{m2}
\\
\vdots & \vdots & \ddots & \vdots
\\
a_{1n} & a_{2n} & ... & a_{mn}
\end{array}
\right\}
\label{eq3:1}
\end{equation}

\noindent
where the first index $(k)$ of an element $a_{kl},\,k=1,2,...,m,$
denotes a column and the second index $(l),\, l=1,2,...,n$,
denotes a row, and where wavy brackets are used
to distinguish a cracovian table of elements from a matrix table of the same elements,
for which ordinary brackets are used.
\\
\noindent
{\bf Definition 3.2.} The product of two cracovians $A$ and $B$, denoted by 
$A\cdot B$ is obtained by multiplying the columns of $A$ by the columns of $B$.
This type of product can be performed only when the 
two cracovians have the same number of rows. The element in the $k-$th column and in the
$l-$th row of the cracovian $A\cdot B$ is obtained by multiplying the $k-$th column of 
$A$ by the $l-$th column of $B$, hence

\begin{equation}
(A\cdot B)_{kl}=\sum\limits_{i} a_{ki}\,b_{li}
\label{eq3:2}
\end{equation}

\noindent
where the summation over the index $i$ extends over all rows.
\\
\noindent
{\bf Observation 3.2.} It will be shown that the cracovian product in Eq. (3.2) leads to the definition of cracovian quasigroups, having the properties of the quasigroups defined in  Section 2. This justifies the use of the dot in the definition in Eq. (3.2). 
\\
\noindent
{\bf Observation 3.3.} The cracovian product can also be defined in the terms
of row-by-row multiplication \cite{Kocinski3}.
\\
\noindent
{\bf Corollary 3.1.} From Eq. (3.2) it follows that the multiplication of cracovians is noncommutative.
\\
\noindent
{\bf Definition 3.3.} The square cracovian $T$ which is called the transposing
cracovian is defined by its elements $t_{kj}$,

\begin{equation}
t_{kj}=1, \quad\mbox{if}\quad k=j,\quad \mbox{and}\quad t_{kj}=0, \quad
\mbox{if}\quad k\neq j, \quad k,j=1,...,n
\label{eq3:3}
\end{equation}

\noindent
The number of rows $n$ of the square cracovian $T$, is equal to the number of rows
of the cracovian which is multiplied by $T$ from the left or from the right. 
\\
\noindent
{\bf Corollary 3.2}. Any cracovian $A$ multiplied by the transposing cracovian $T$ from the right remains
unchanged, and it is changed to the transposed cracovian after the multiplication 
by $T$ from the left, since we have:

\begin{equation}
(A\cdot T)_{kl}=\sum\limits_i a_{ki}t_{li}=a_{kl}
\label{eq3:4}
\end{equation}

\begin{equation}
(T\cdot A)_{kl}=\sum\limits_i t_{ki}a_{li}=a_{lk}
\label{eq3:5}
\end{equation}

\noindent
{\bf Corollary 3.3}. The transposing cracovian $T$ is the right unit cracovian, however, 
it is not a left unit cracovian.
\\ 
\noindent
{\bf Corollary 3.4.}  We have:
\begin{equation}
T\cdot T=T
\label{eq3:6}
\end{equation}

\noindent
{\bf Corollary 3.5.} The cracovian $T$ is the only cracovian having the 
properties specified in Eqs. (\ref{eq3:4}), (\ref{eq3:5}) and (\ref{eq3:6}).
\\
\noindent
{\bf Definition 3.5.} The multiplication of cracovians is performed from left to right:

\begin{equation}
A\cdot B\cdot C\cdot D \ldots\cdot Z:=\{[(A\cdot B)\cdot C]\cdot D\cdot\ldots\cdot Z\}
\label{eq3:7}
\end{equation}

\noindent
This multiplication is possible when the number of rows of a particular factor in 
the chain of factors is equal to
the number of rows of the preceding resultant cracovian.
\\
\noindent 
{\bf Corollary 3.6.} For any cracovian $A$ we have the equality: 
\begin{equation}
T\cdot (T\cdot A)=A
\label{eq3:8}
\end{equation}

\noindent
{\bf Corollary 3.7.} For any two cracovians $A$ and $B$ with the same number of rows
we have the equality:

\begin{equation}
T\cdot (A\cdot B)=(B\cdot A)
\label{eq3:9}
\end{equation}

\noindent
Proof.
According to Eq. (3.5), the element $(kl)$ of $A\cdot B$ turns into 
the element $(lk)$ of $T\cdot (A\cdot B)$, and the latter is equal to the element
$(lk)$ of $B\cdot A$.
\\
\noindent
{\bf Corollary 3.8.} The square $A^2=A\cdot A$, of any nonzero cracovian
$A$ is a symmetric cracovian, since

\begin{equation}
T\cdot (A\cdot A)=A\cdot A
\label{eq3:10}
\end{equation}

\noindent
{\bf Definition 3.5.} In order to avoid unnecessary brackets in the formulas, from now on
we will omit the dot $(\cdot)$ in the product of any cracovian with the transposing
cracovian $T$, from the left side or from the right side. This means that we will write
$TA$ instead of $T\cdot A$, and $AT$ instead of $A\cdot T$. Further on we will write 
$B\cdot TA$ instead of $B\cdot (T\cdot A)$, and $AT\cdot B$ will replace
$(A\cdot T)\cdot B$. 
\\
\noindent
{\bf Lemma 3.1.} \cite{Kocinski3}.
For three cracovians $A, B$ and $C$ we have the equality:

\begin{equation}
(A\cdot B)\cdot C=A\cdot (C\cdot TB)
\label{eq3:11}
\end{equation}

\noindent 
{\bf Corollary 3.9}. From Eqs. (\ref{eq3:9}) and (\ref{eq3:11}), it follows that

\begin{equation}
T[(A\cdot B)\cdot C]=(C\cdot TB)\cdot A
\label{eq3:12}
\end{equation}

\noindent
{\bf Corollary 3.10}. For an arbitrary finite number of factors we obtain,

\begin{equation}
T(A_1\cdot A_2\cdot ... \cdot A_{k-1}\cdot A_k)=A_k\cdot TA_{k-1}\cdot ...\cdot TA_2\cdot A_1
\label{eq3:13}
\end{equation}

\noindent
{\bf Corollary 3.11}. For the column-by-column product there is no 
cracovian with the property of the left unit cracovian.
\\
\noindent
Proof. If $X$ were the left unit cracovian, we would have: $X\cdot (B\cdot C)
=B\cdot C$. From Eq. (\ref{eq3:11}) at the same time we obtain: $X\cdot (B\cdot C)=
[X\cdot (TC)]\cdot B=TC\cdot B$. We therefore obtain the equality: $B\cdot C=TC\cdot B$,
which in a general case is not true.
\\
\noindent
{\bf Corollary 3.12}. The relations between the matrix product and the cracovian product of two tables $A$ and $B$, and of three tables $A$, $B$ and $C$, are given by the respective equalities:

\begin{equation}
AB=B\cdot TA,\qquad {\rm and}\qquad ABC=C\cdot (TB)\cdot (TA)
\label{eq3:14}
\end{equation}

\begin{equation}
{\tilde B}A=A\cdot B,\qquad {\rm and}\qquad \tilde C{\tilde B}A=A\cdot B\cdot C
\label{eq3:15}
\end{equation}

\noindent
where on the left hand side of these equalities, the tables undergo the matrix product, 
and on the right hand side the same tables undergo the cracovian product, 
and where the symbol $\tilde{\null}$ denotes the transposed matrix \cite{Kocinski1,Kocinski3}.
\\
\noindent
{\bf Definition 3.6}. The right inverse of a square cracovian $A$ is defined as the cracovian
$A^{-1}$ for which we have:

\begin{equation}
A\cdot A^{-1}=T
\label{eq3:15}
\end{equation}

\noindent
From Eqs. (3.8) and (3.11) it follows that $A^{-1}$ is also the left inverse.
\\
\noindent
{\bf Corollary 3.13.} For a square cracovian $A$ with the inverse $A^{-1}$, we have

\begin{equation}
(TA)\cdot (TA^{-1})=T
\label{eq3:16}
\end{equation}

\noindent
{\bf Corollary 3.14.} For the product of two square cracovians $A_1$ and $A_2$, having the inverses we have:

\begin{equation}
(A_1\cdot A_2)^{-1}=A_1^{-1}\cdot A_2^{-1}
\label{eq3:17}
\end{equation}

\noindent
{\bf Corollary 3.15.} From Eq. (3.18) we obtain:

\begin{equation}
(TA)^{-1}=TA^{-1}
\label{eq3:18}
\end{equation}

\noindent
{\bf Observation 3.4}. A square table of symbols or complex numbers has the cracovian inverse if and only if it has the matrix inverse.
\\
\noindent
{\bf Lemma 3.2} \cite{Kocinski3}. Between the matrix inverse $A^{-1}_m$ of a square table of symbols or of complex numbers, and the cracovian inverse $A^{-1}_c$ of that table holds the relation: 

\begin{equation}
A^{-1}_m=T(A^{-1}_c)
\label{eq3:19}
\end{equation}

\noindent
The matrix inverse of a square table is equal to the transposed cracovian inverse of that table.

\vspace{5mm}
\noindent
In the following formulas
we start with the matrix expressions, and employing Eqs. (3.15) and (3.20), we determine the respective cracovian expressions.
\\
\noindent
{\bf Corollary 3.16}. Let $\vec e_m$ and $\vec e_c$ denote the one-column matrix and the one-column  cracovian, respectively, constructed from the basis vectors $\vec e_1,...,\vec e_n$, of an $n-$dimensional linear vector space. Let $S_m$ and $S_c$ be the respective matrix and cracovian transformation of these basis vectors. The change of basis then is defined by the expressions:

\begin{equation}
\vec e^{\,\prime}_m=\tilde S_m\vec e_m=\vec e_c\cdot T(TS_c)=\vec e_c\cdot S_c=\vec e^{\,\prime}_c
\label{eq3:21}
\end{equation}

\noindent
where $\tilde S_m$ denotes the transposed matrix, and where in the first equality we are dealing with the matrix product, while in the next two equalities we are dealing with the cracovian product.
\\
\noindent
{\bf Corollary 3.17}. Let $x_m$ and $x_c$ (or $y_m$ and $y_c$) denote the one-column matrix (or the one-column cracovian) constructed from the components $x_1,x_2,...,x_n$ (or $y_1,y_2,...,y_n$), of a vector, which is referred to the two bases in Eq. (3.21), respectively. The relation between the two sets of components is given by

\begin{equation}
x_m=S_my_m=y_c\cdot TS_c=x_c
\label{eq3:22}
\end{equation}

\noindent
where the transformation table $S$, has been identified  with the matrix $S_m$ or with the cracovian $S_c$, depending on the type of the employed product.
\\
\noindent
{\bf Corollary 3.18}. Let $x_{\rm m}$ and $x_{\rm c}$ denote the column matrix and column
cracovian, respectively, constructed from the components
$(x_1, x_2,...,x_n)$ of a vector in an $n-$dimensional linear vector space.
A linear mapping of that vector, expressed in the terms of matrices or cracovians 
is given by

\begin{equation}
x^{\prime}_{\rm m}=A_{\rm m}\,x_{\rm m}=x_{\rm c}\cdot TA_{\rm c}=x^{\prime}_{\rm c}
\label{eq3:23}
\end{equation}

\noindent
where the table $A$ is identified with a matrix $A_{\rm m}$ or with a cracovian
$A_{\rm c}$, depending on the type of the employed product. 
\\
\noindent
{\bf Corollary 3.19.} The relation between matrix and cracovian linear transformations
$A_{\rm m}$ and $B_{\rm m}$, or $A_{\rm c}$ and $B_{\rm c}$, respectively, which determine
the same linear mapping referred to two bases which are 
connected according to Eq. (3.21) is given by

\begin{equation}
B_{\rm m}=S^{-1}_{\rm m}\,A_{\rm m}\,S_{\rm m}=S_{\rm c}\,\cdot TA_{\rm c}\,
\cdot S^{-1}_{\rm c}=B_{\rm c}
\label{eq3:24}
\end{equation}

\noindent
where the second equality connects an expression in matrix product with the same expression in cracovian product, and that connection follows from the second of Eqs. (3.15) and from Eq. (3.20).
\\
\noindent
{\bf Corollary 3.20}. From the commutation condition of two matrix tables $A_m$ and $B_m$:
$A_mB_m=B_mA_m$, and the relation between matrix product and cracovian product of two tables
$A$ and $B$ inthe first of equations in Eq. (3.15), we obtain the commutativity condition of two cracovians $A_c$ and $B_c$ in the form:

\begin{equation}
A_c\cdot T B_c=B\cdot TA_c
\label{eq3:25}
\end{equation}

\section{Cracovian representations of quasigroups connected with groups with involution}

\noindent
{\bf Observation 4.1}. To the definition of the "dot" product in Eq. (2.3)
coresponds the relation between cracovian and matrix product in the first of the two relations in Eq. (3.15).
The consequences following from the product in Eq. (2.3) are deduced in Eqs. (2.4) through (2.20). The consequences which follow from the cracovian product defined in Eq. (3.2)
are deduced in Eqs. (3.3) through (3.20); the first of the two relations between cracovian product and matrix product in Eq. (3.15), appears among them. To each consequence of the "dot"
product in Eq. (2.3) there corresponds the respective consequence from the cracovian product in Eq. (3.2).
We conclude that there exists a homomorphism between the 
quasigroups defined by the product in Eq. (2.3), and square cracovians with 
inverses, provided that the involution operation $I$ in Eq. (2.2) is identified with the operation of transposing a matrix. Consequently, square
cracovians with inverses constitute representations of quasigroups connected with
groups with involution. The right identity $\tau$ is represented by the transposing cracovian $T$.
\\
\noindent
{\bf Corollary 4.1}. To	every matrix group with the involution operation, which is 
the operation of transposing a matrix, corresponds a cracovian quasigroup.
\\
\noindent
{\bf Definition 4.1.}  A representation by linear substitutions of a quasigroup $QG$ 
connected with a group $G$ in Eq. (2.1), is a cracovian quasigroup
onto which the quasigroup is homomorphic. It consists of the assignment of a square
cracovian $C(\hat a)$ to each quasigroup element $\hat a$ in such a way that,

\begin{equation}
C\hat a)\cdot C(\hat b)=C(\hat a\cdot\hat b), \qquad\forall \hat a,\hat b\in QG
\label{eq4:1}
\end{equation}

\noindent
{\bf Observation 4.2}. The expression "a reducible set of matrices" is used in the sense of "a completely reducible set of matrices", i. e. to a set of matrices which is equivalent to the direct sum of two or more other sets of matrices. The expression "an irreducible set of matrices" is used in the sense of a set of matrices which is not equivalent to the direct sum of two or more other sets of matrices.
\\
\noindent
{\bf Observation 4.3.} The notion of irreducibility of a matrix representation of a group
carries over on a cracovian representation of a quasigroup. However, if
a set of square tables of complex numbers in its quality of
being a cracovian  irrep of a quasigroup $QG$ in Eq. (4.1),
were identified with the matrices representing the group with involution $G$, with which 
the respective quasigroup $QG$ is connected, that set of matrices could be reducible. 
\cite{Kocinski1,Kocinski2}.
\\
\noindent
{\bf Observation 4.4}. The notions of faithful or unfaithful matrix representations
defined for groups, carry over on cracovian representations of quasigroups connected
with groups with involution. Each cracovian quasigroup is its own faithful representation.
\\
\noindent
{\bf Observation 4.5}. Two types of cracovian irreps of quasigroups 
connected with groups with involution are known. (1) There are cracovian
irreps of quasigroups connected with groups with involution,
which at the same time are matrix irreps of the respective
groups with involution. To these belong the single-valued cracovian irreps of the quasigroups connected with the orthogonal and pseudo-orthogonal continuous rotation groups.
(2) There are double-valued cracovian irreps of quasigroups connected with groups with involution 
which are not matrix irreps of those groups at the same time. 
To these belong the double-valued cracovian irreps of the orthogonal and pseudo-orthogonal rotation groups $SO(3),\,L^{\uparrow}_+,\,SO(3,2),\,SO(4,1)$
\cite{Kocinski1,Kocinski2,Kocinski3}.

\section
{\bf The "hidden" symmetry group of Weyl, and the quasigroup connected with a group with involution}

\noindent
It seems that an analogy can be drawn between the "hidden" symmetry group of an object, introduced by Weyl \cite{Weyl}, and the quasigroup, connected with a group with involution
\cite{Kocinski2}.
\\
\indent
Weyl considered a set with a symmetry group $G$. This can be the set of all roots of a polynom, the set of space-time points,  or of all nodes of a crystal lattice. He showed that the essential features of a set endowed with a structure can be determined by studying the group of all automorphisms $Aut\,G$ of this set, which preserve all structural relations. The group $G$ determines the "obvious" symmetry, and the group $Aut\,G$, the "hidden" symmetry of the set. The concept of a "hidden" symmetry group was discussed from the standpoint of group actions on sets by Florek {\em et al.}\cite{Florek}, and by Lulek {\em et al.}\cite{Lulek}.
\\
\indent
In one of the examples discussed by Weyl, the object is a regular septadecagon. The property under consideration is the possibility of construction of that regular septadecagon with the help of a compass and a ruler. The symmetry group $G$ of the set of vertices of the regular septadecagon is the cyclic group $C_{17}$. This is the obvious geometric symmetry group of the regular septadecagon. The vertices of a regular septadecagon are determined by the roots of the equation $z^{17}-1=0$, with $z=x+iy$. One of the roots is $z=1$, and the remaining 16 roots are determined by an algebraic equation of degree 16. The root $z=1$ determines the starting vertex in the construction of a regular septadecagon. The determination of the positions of the remainig 16 vertices is shown to be connected with the group of permutations of the 16 roots of the algebraic equation of the degree 16. This appears to be the cyclic group $C_{16}$. Consequently, $Aut\,C_{17}=C_{16}$, which is the hidden symmetry group of a regular septadecagon. The possibility of construction of a regular septadecagon with a compass and a ruler hinges on the group $C_{16}$.
\\
\indent 
The group of automorphisms $I(g)$ of the obvious symmetry group $G$ determines the quasigroup $QG$.
The quasigroup $QG$ connected with a group with involution $G$, could be considered as a hidden symmetry group $Aut\,G$ of the object whose obvious symmetry is determined by the 
group G.
The group of automorphisms defined in Eq. (2.2) is the basis for defining the nonassociative product in Eq.(2.3). This product serves for the definition of a quasigroup which could
be recognized  as a particular case of Weyl's hidden symmetry group of an object.

\section{\bf Conclusions}

\noindent
It has been shown that in groups with involution an automorphism can be defined, which leads to the definition of a nonassociative product. The group elements undergoing this nonassociative product form a certain type of a quasigroup, which has the right unit element but not a left unit element.
If in the matrix representation of a group with involution, the operation of involution is identified with the operation of transposition of a matrix,
the quasigroup connected with that group has cracovian representations. A cracovian representation is analogous to a matrix representation with two basic differences: The row-by-column product of two matrices is replaced by the column-by-column product of two cracovians, and there is only the right unit element. The orthogonal and pseudo-orthogonal rotation groups are examples of groups with involution. In the matrix representations of these groups any transposed matrix belongs to the relevant group. The involution operation in the matrix group can therefore be identified with the transposition of a matrix, and we can define the nonassociative product, which can be identified with the cracovian product.
An analogy has been drawn between Weyl's
"hidden" symmetry group and a quasigroup.

\end{document}